\documentclass[12pt,leqno,draft]{article}

\begin{document}

\title{Analysis on disconnected sets}

\author{Stephen Semmes \\
        Rice University}

\date{}

\maketitle

        Very often in analysis, one focuses on connected spaces.  This
is certainly not always the case, and in particular there are many
interesting matters related to Cantor sets.  Here we are more
concerned with a type of complementary situation.  As a basic
scenario, suppose that $U$ is an open set in ${\bf R}^n$ and that $E$
is a closed set contained in the boundary of $U$ such that for every
$x \in U$ and $r > 0$ there is a connected component of $U$ contained
in the ball with center $x$ and radius $r$.  For instance, $E$ might
be the boundary of $U$.  As a uniform version of the condition, one
might ask that there be a constant $C > 0$ such that the
aforementioned connected component of $U$ contains a ball of radius
$C^{-1} \, r$.  The connected components of $U$ might be quite
regular, even if there are infinitely many of them.

        As a basic example, $E$ could be a Cantor set in the real
line, and $U$ could be the complement of $E$ or the union of the
bounded complementary components of $E$.  This does not work in higher
dimensions, where the complement of a Cantor set is connected.
Sierpinski gaskets and carpets in the plane are very interesting cases
where $E$ is connected.  One could also consider more abstract
versions of this, in metric spaces, for instance.  The open set $U$
could be discrete, and the limiting set $E$ could be arbitrary.

        Let us restrict our attention to the case where the closure of
$U$ is compact and $E = \partial U$.  Consider Bergman spaces of
locally constant functions on $U$, which can be identified with
$\ell^2$ spaces.  Toeplitz operators associated to continuous
functions on $\overline{U}$ can be defined in the usual way, by
multiplying and projecting.  The Toeplitz operators associated to
functions vanishing on the boundary are compact.  Thus, modulo compact
operators, it is really the functions on the boundary that are
important.  Invertible functions on the boundary correspond to
Fredholm operators, but their indices are automatically equal to $0$.
Of course, these are well-known themes.

        There are fancier versions of this using Bergman spaces of
holomorphic functions.  One could also look at analogous notions based
on Hardy spaces, with boundary norms on the individual components.
Note that there are two types of holomorphic functions in Clifford
analysis, because of noncommutativity, and the product of a
holomorphic function and a constant is holomorphic when the product is
in the appropriate order.  When multiplying a Clifford holomorphic
function by a continuous function, it is better to have the continuous
function in the same place as the constant in the previous statement.

        Analogous to the disk algebra, there is the space of
continuous functions on $\overline{U}$ which are holomorphic on $U$.
This can be quite nontrivial for regions in the plane with suitably
fractal boundary.  By contrast, a continuous function on a
neighborhood of a nice arc in the plane which is holomorphic on both
sides of the arc is holomorphic across the arc.  Fractal sets can also
be removed when the function is sufficiently regular.  The product of
a continuous function on $\overline{U}$ which is holomorphic on $U$
with a continuous function on $\overline{U}$ is continuous on
$\overline{U}$ but not normally holomorphic on $U$.  However, the
product can be corrected to give a continuous function on
$\overline{U}$ which is holomorphic on $U$ by solving a
$\overline{\partial}$ probem when the second function is smooth
enough.  The relationship among the components of $U$ through
continuity seems intriguing.

        There certainly seem to be some issues here that are not
quite the same as in the classical situation.

        The components of $U$ are related to topological activity in
the boundary, and infinitely many components indicate a lot of
activity.  It is natural to look at analysis that reflects this.

\end{document}